\documentclass[a4paper,12pt]{amsart}
\usepackage{tabularx}
\usepackage[dvips]{graphicx, color}
\usepackage{amsmath}
\usepackage{amssymb}
\usepackage{amsbsy}

\usepackage{amsgen}
\usepackage{amsopn}
\usepackage{amscd}
\usepackage{amstext}
\usepackage{amsthm}

\theoremstyle{definition}

\theoremstyle{plain}

\theoremstyle{remark}

\theoremstyle{definition}

\newtheorem{acknow}{Acknowledgment}
\renewcommand{\acknow}{\textbf{Acknowledgments--}}
\newtheorem{conv}{Conventions}
\renewcommand{\conv}{\textbf{Conventions}.  }

\begin{document}
\title{On Reidemeister invariance of the Khovanov homology group of the Jones polynomial}
\author{Noboru Ito}
\maketitle

\begin{abstract}As Oleg Viro describes in his paper, the most fundamental property of the Khovanov homology group is their invariance under Reidemeister moves.  Viro constructes Khovanov complex and homology consisting of Jordan curves with sign and also gives a proof for the only case of first Reidemeister move by using his definition of Khovanov homology groups.  In this paper, homotopy maps are obtained explicitly for the other Reidemeister moves, i.e. second and third.  
\end{abstract}

\section{Introduction.  }
\subsection{Motivation.  }
Oleg Viro makes a construction of the Khovanov homology group of Jones polynomial significantly simpler \cite{viro}.  The chain maps induced by a Reidemeister move are homotopy equivalences.  Thus, there exist chain homotopy maps that deduces the homotopy equivalences.  Magnus Jacobsson obtains explicitly chain maps induced by a link cobordism \cite{jacobsson}.  
However, these homotopy maps are missing.  Therefore, Viro gives a proof for the only case of first Reidemeister move \cite{viro}.  
In this paper, homotopy mappings are obtained explicitly for the other Reidemeister moves, i.e. second and third, by using Viro's definition of Khovanov homology groups and chain maps \cite{jacobsson}.

\subsection{Short review of Viro's definition of the Khovanov homology group using Jacobsson's description.  }

A version of Jones polynomial is defined by (1)--(3).  
\begin{align}
&V_{L}(unknot) = q + q^{-1}, \\
&q^{-2}V_{L}\left(\begin{minipage}{15pt}
\begin{picture}(15,15)
\put(0,0){\vector(1,1){15}}
\qbezier(15,0)(15,0)(10,5)
\qbezier(5,10)(0,15)(0,15)
\put(1,14){\vector(-1,1){1}}
\end{picture}
\end{minipage}
\right) - q^{2}V_{L}\left(\begin{minipage}{15pt}
\begin{picture}(15,15)
\put(15,0){\vector(-1,1){15}}
\qbezier(0,0)(5,5)(5,5)
\qbezier(10,10)(15,15)(15,15)
\put(14,14){\vector(1,1){1}}
\end{picture}
\end{minipage}
\right) = (q^{-1} - q) V_{L}\left(
\begin{minipage}{15pt}
\begin{picture}(15,15)
\qbezier(0,0)(10,7.5)(0,15)
\qbezier(15,0)(5,7.5)(15,15)
\put(14,14){\vector(1,1){1}}
\put(1,14){\vector(-1,1){1}}
\end{picture}
\end{minipage}\right).  
\end{align}

By using the definitions, the version of Jones polynimilal is represented in the following manner as \cite[Section 2]{jacobsson} or \cite[Section 4]{viro}.  Let $D$ be diagram of Link $L$.  

\begin{equation}
V_{L} = \sum_{\text{Kauffman states $S$ of} D} (-1)^{\frac{w(D)-\sigma(S)}{2}}q^{\frac{3w(D) - \sigma(S)}{2}}(q + q^{-1})^{r_{s}}
\end{equation}
where $r_{s}$ is the number of states.  

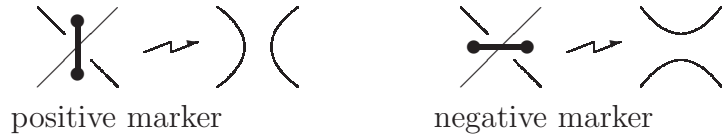
\begin{figure}[!hbtp]
\begin{center}
\begin{minipage}{300pt}
\begin{picture}(40,40)
\put(0,10){\line(1,1){30}}
\qbezier(30,10)(30,10)(20,20)
\qbezier(10,30)(0,40)(0,40)
\qbezier(40,23)(45,25)(50,27)
\qbezier(50,27)(50,25)(50,23)
\qbezier(50,23)(55,25)(60,27)
\put(60,27){\vector(4,1){0}}
\linethickness{2pt}
\put(15,15){\circle*{5}}
\put(15,35){\circle*{5}}
\put(15,15){\line(0,1){20}}
\put(-10, -5){\text{positive marker}}
\end{picture}
\qquad
\begin{picture}(40,40)
\qbezier(0,10)(20,25)(0,40)
\qbezier(30,10)(10,25)(30,40)
\end{picture}
\qquad\qquad
\begin{picture}(40,40)
\put(0,10){\line(1,1){30}}
\qbezier(30,10)(30,10)(20,20)
\qbezier(10,30)(0,40)(0,40)
\qbezier(40,23)(45,25)(50,27)
\qbezier(50,27)(50,25)(50,23)
\qbezier(50,23)(55,25)(60,27)
\put(60,27){\vector(4,1){0}}
\linethickness{2pt}
\put(5,25){\circle*{5}}
\put(25,25){\circle*{5}}
\put(5,25){\line(1,0){20}}
\put(-10, -5){\text{negative marker}}
\end{picture}
\qquad
\begin{picture}(40,40)
\qbezier(0,10)(15,30)(30,10)
\qbezier(0,40)(15,20)(30,40)
\end{picture}
\end{minipage}
\caption{Smoothing of a diagram according to thick segments corresponding to markers.  }
\label{marker1}
\end{center}
\end{figure}

For the convenience, we use the terminology and notation in \cite[Page 1213, Subection 2.1]{jacobsson}, which are a {\it marker} (See Figure \ref{marker1}), {\it ``refined" states}, $\tau(S)$, $i(S)$ and $j(S)$.  By using notation, we have
\begin{equation}
V_{L} = \sum_{\text{``refined" states $S$ of} D} (-1)^{i(S)}q^{j(S)}.  
\end{equation}

In the rest of this notes, we use the terminology, notation or symbol from \cite[Page 1213--1218, Subection 2.1--2.3]{jacobsson}.   

\section{Proof of the Reidemeister invariance.  }
The proof for the case of first Reidemeister move is obtained by O. Viro \cite{viro}.  Then, we prove for the other moves in the manner as in \cite{viro} using the notation in \cite[Figure 9, 10, 13, 14, 17, 18 and Section 3.3.1]{jacobsson}.  
\subsection{Second Reidemeister move}
Let $a$, $b$ be crossings, $x$ sequence of crossings with negative markers  and $p$, $q$ be signs.  For a crossing with no markers or no signs in the following formulas, any markers or signs may be selected.  We use the symbols $p:q$ and $q:p$ whose definitions are given by \cite[Figure 3 Subsection 2.2]{jacobsson}.  
\begin{align}\label{2nd-sum}
& \mathcal{C}\left(~
    \begin{minipage}{30pt}
        \begin{picture}(30,30)
            \qbezier(0,0)(40,15)(0,30)
            \qbezier(30,30)(8,20)(30,30)
            \qbezier(12,20)(2,15)(12,10)
            \qbezier(30,0)(10,10)(30,0)
            \put(4,20){$a$}
            \put(4,5){$b$}
        \end{picture}
    \end{minipage}
~\right) = \\
& 
\mathcal{C}\left(~\begin{minipage}{40pt}
        \begin{picture}(30,40)
            \qbezier(0,5)(40,20)(0,35)
            \qbezier(30,35)(8,25)(30,35)
            \qbezier(11,26)(1,20)(11,14)
            \qbezier(30,5)(10,15)(30,5)
  {\color{red}{\put(14,33){\circle*{3}}
\put(14,24){\circle*{3}}
\put(14,24){\line(0,1){8}}}}
            {\color{blue}{\put(10,17){\circle*{3}}
\put(10,8){\circle*{3}}
\put(10,9){\line(0,1){8}}}}
\put(-5,19){$p$}
\put(18,19){$q$}
        \end{picture}
    \end{minipage} \otimes [xa] + \begin{minipage}{40pt}
        \begin{picture}(30,40)
            \qbezier(0,5)(40,20)(0,35)
            \qbezier(30,35)(8,25)(30,35)
            \qbezier(11,26)(1,20)(11,14)
            \qbezier(30,5)(7,14)(30,5)
            {\color{blue}
            {\put(8,27.5){\circle*{3}}
\put(20,27.5){\circle*{3}}
\put(9,27.5){\line(1,0){10}}
            }}
            {\color{red}{\put(16,12){\circle*{3}}
\put(5,12){\circle*{3}}
\put(5,12){\line(1,0){10}}}}
\put(5,17){\text{$-$}}
\put(1,37){$p:q$}
\put(0,0){$q:p$}
        \end{picture}
    \end{minipage} \otimes [xb]~\right)
    \oplus
    \mathcal{C}\left(~ 
    \begin{minipage}{30pt}
        \begin{picture}(30,30)
            \qbezier(0,0)(40,15)(0,30)
            \qbezier(30,30)(8,20)(30,30)
            \qbezier(11,21)(1,15)(11,9)
            \qbezier(30,0)(10,10)(30,0)
            {\color{blue}
            {\put(8,22.5){\circle*{3}}
\put(20,22.5){\circle*{3}}
\put(9,22.5){\line(1,0){10}}
            }}
            {\color{blue}{\put(10,12){\circle*{3}}
\put(10,3){\circle*{3}}
\put(10,4){\line(0,1){8}}}}
        \end{picture}
    \end{minipage} \otimes [x], \begin{minipage}{30pt}
        \begin{picture}(30,30)
            \qbezier(0,0)(40,15)(0,30)
            \qbezier(30,30)(8,20)(30,30)
            \qbezier(11,21)(1,15)(11,9)
            \qbezier(30,0)(7,9)(30,0)
  {\color{red}{\put(14,28){\circle*{3}}
\put(14,19){\circle*{3}}
\put(14,19){\line(0,1){8}}}}
            {\color{red}{\put(16,7){\circle*{3}}
\put(5,7){\circle*{3}}
\put(5,7){\line(1,0){10}}}}
        \end{picture}
    \end{minipage} \otimes [xab] , 
    \begin{minipage}{30pt}
        \begin{picture}(30,30)
            \qbezier(0,0)(40,15)(0,30)
            \qbezier(30,30)(8,20)(30,30)
            \qbezier(11,21)(1,15)(11,9)
            \qbezier(30,0)(7,9)(30,0)
            {\color{blue}
            {\put(8,22.5){\circle*{3}}
\put(20,22.5){\circle*{3}}
\put(9,22.5){\line(1,0){10}}
            }}
            {\color{red}{\put(16,7){\circle*{3}}
\put(5,7){\circle*{3}}
\put(5,7){\line(1,0){10}}}}
\put(5,12){\text{$-$}}
        \end{picture}
    \end{minipage} \otimes [xb] ~\right).  \nonumber
\end{align}

The isomorphism 

$\mathcal{C}\left(~\begin{minipage}{40pt}
        \begin{picture}(30,40)
            \qbezier(0,5)(40,20)(0,35)
            \qbezier(30,35)(8,25)(30,35)
            \qbezier(11,26)(1,20)(11,14)
            \qbezier(30,5)(10,15)(30,5)
  {\color{red}{\put(14,33){\circle*{3}}
\put(14,24){\circle*{3}}
\put(14,24){\line(0,1){8}}}}
            {\color{blue}{\put(10,17){\circle*{3}}
\put(10,8){\circle*{3}}
\put(10,9){\line(0,1){8}}}}
\put(-5,19){$p$}
\put(18,19){$q$}
        \end{picture}
    \end{minipage} \otimes [xa] + \begin{minipage}{40pt}
        \begin{picture}(30,40)
            \qbezier(0,5)(40,20)(0,35)
            \qbezier(30,35)(8,25)(30,35)
            \qbezier(11,26)(1,20)(11,14)
            \qbezier(30,5)(7,14)(30,5)
            {\color{blue}
            {\put(8,27.5){\circle*{3}}
\put(20,27.5){\circle*{3}}
\put(9,27.5){\line(1,0){10}}
            }}
            {\color{red}{\put(16,12){\circle*{3}}
\put(5,12){\circle*{3}}
\put(5,12){\line(1,0){10}}}}
\put(5,17){\text{$-$}}
\put(1,37){$p:q$}
\put(0,0){$q:p$}
        \end{picture}
    \end{minipage} \otimes [xb]~\right)$
$\to$ $\mathcal{C} \left(~
    \begin{minipage}{30pt}
        \begin{picture}(30,30)
            \qbezier(0,0)(20,15)(0,30)
            \qbezier(30,0)(10,15)(30,30)
        \end{picture}
    \end{minipage}\otimes [x]
~\right)$
\vspace{3mm}

is defined by the formulas

\begin{equation}
\begin{minipage}{40pt}
        \begin{picture}(30,40)
            \qbezier(0,5)(40,20)(0,35)
            \qbezier(30,35)(8,25)(30,35)
            \qbezier(11,26)(1,20)(11,14)
            \qbezier(30,5)(10,15)(30,5)
  {\color{red}{\put(14,33){\circle*{3}}
\put(14,24){\circle*{3}}
\put(14,24){\line(0,1){8}}}}
            {\color{blue}{\put(10,17){\circle*{3}}
\put(10,8){\circle*{3}}
\put(10,9){\line(0,1){8}}}}
\put(-5,19){$p$}
\put(18,19){$q$}
        \end{picture}
    \end{minipage} \otimes [xa] + \begin{minipage}{40pt}
        \begin{picture}(30,40)
            \qbezier(0,5)(40,20)(0,35)
            \qbezier(30,35)(8,25)(30,35)
            \qbezier(11,26)(1,20)(11,14)
            \qbezier(30,5)(7,14)(30,5)
            {\color{blue}
            {\put(8,27.5){\circle*{3}}
\put(20,27.5){\circle*{3}}
\put(9,27.5){\line(1,0){10}}
            }}
            {\color{red}{\put(16,12){\circle*{3}}
\put(5,12){\circle*{3}}
\put(5,12){\line(1,0){10}}}}
\put(5,17){\text{$-$}}
\put(1,37){$p:q$}
\put(0,0){$q:p$}
        \end{picture}
    \end{minipage} \otimes [xb]
\mapsto 
    \begin{minipage}{30pt}
        \begin{picture}(30,30)
            \qbezier(0,0)(20,15)(0,30)
            \qbezier(30,0)(10,15)(30,30)
            \put(0,15){$p$}
\put(23,15){$q$}
        \end{picture}
    \end{minipage}\otimes [x].  
\end{equation}

The retraction $\rho : \mathcal{C}
\left(~
    \begin{minipage}{30pt}
        \begin{picture}(30,30)
            \qbezier(0,0)(40,15)(0,30)
            \qbezier(30,30)(8,20)(30,30)
            \qbezier(12,20)(2,15)(12,10)
            \qbezier(30,0)(10,10)(30,0)
            \put(4,20){$a$}
            \put(4,5){$b$}
        \end{picture}
    \end{minipage}
~\right)$ $\to$
$\mathcal{C}\left(~\begin{minipage}{40pt}
        \begin{picture}(30,40)
            \qbezier(0,5)(40,20)(0,35)
            \qbezier(30,35)(8,25)(30,35)
            \qbezier(11,26)(1,20)(11,14)
            \qbezier(30,5)(10,15)(30,5)
  {\color{red}{\put(14,33){\circle*{3}}
\put(14,24){\circle*{3}}
\put(14,24){\line(0,1){8}}}}
            {\color{blue}{\put(10,17){\circle*{3}}
\put(10,8){\circle*{3}}
\put(10,9){\line(0,1){8}}}}
\put(-5,19){$p$}
\put(18,19){$q$}
        \end{picture}
    \end{minipage} \otimes [xa] + \begin{minipage}{40pt}
        \begin{picture}(30,40)
            \qbezier(0,5)(40,20)(0,35)
            \qbezier(30,35)(8,25)(30,35)
            \qbezier(11,26)(1,20)(11,14)
            \qbezier(30,5)(7,14)(30,5)
            {\color{blue}
            {\put(8,27.5){\circle*{3}}
\put(20,27.5){\circle*{3}}
\put(9,27.5){\line(1,0){10}}
            }}
            {\color{red}{\put(16,12){\circle*{3}}
\put(5,12){\circle*{3}}
\put(5,12){\line(1,0){10}}}}
\put(5,17){\text{$-$}}
\put(1,37){$p:q$}
\put(0,0){$q:p$}
        \end{picture}
    \end{minipage} \otimes [xb]~\right)$ is defined by the formulas 
    \[\begin{minipage}{40pt}
        \begin{picture}(30,40)
            \qbezier(0,5)(40,20)(0,35)
            \qbezier(30,35)(8,25)(30,35)
            \qbezier(11,26)(1,20)(11,14)
            \qbezier(30,5)(10,15)(30,5)
  {\color{red}{\put(14,33){\circle*{3}}
\put(14,24){\circle*{3}}
\put(14,24){\line(0,1){8}}}}
            {\color{blue}{\put(10,17){\circle*{3}}
\put(10,8){\circle*{3}}
\put(10,9){\line(0,1){8}}}}
\put(-5,19){$p$}
\put(18,19){$q$}
        \end{picture} 
        \end{minipage} \otimes [xa]
        \mapsto \begin{minipage}{40pt}
        \begin{picture}(30,40)
            \qbezier(0,5)(40,20)(0,35)
            \qbezier(30,35)(8,25)(30,35)
            \qbezier(11,26)(1,20)(11,14)
            \qbezier(30,5)(10,15)(30,5)
  {\color{red}{\put(14,33){\circle*{3}}
\put(14,24){\circle*{3}}
\put(14,24){\line(0,1){8}}}}
            {\color{blue}{\put(10,17){\circle*{3}}
\put(10,8){\circle*{3}}
\put(10,9){\line(0,1){8}}}}
\put(-5,19){$p$}
\put(18,19){$q$}
        \end{picture}
    \end{minipage} \otimes [xa] + \begin{minipage}{40pt}
        \begin{picture}(30,40)
            \qbezier(0,5)(40,20)(0,35)
            \qbezier(30,35)(8,25)(30,35)
            \qbezier(11,26)(1,20)(11,14)
            \qbezier(30,5)(7,14)(30,5)
            {\color{blue}
            {\put(8,27.5){\circle*{3}}
\put(20,27.5){\circle*{3}}
\put(9,27.5){\line(1,0){10}}
            }}
            {\color{red}{\put(16,12){\circle*{3}}
\put(5,12){\circle*{3}}
\put(5,12){\line(1,0){10}}}}
\put(5,17){\text{$-$}}
\put(1,37){$p:q$}
\put(0,0){$q:p$}
        \end{picture}
    \end{minipage} \otimes [xb], \]
    \[\begin{minipage}{40pt}
        \begin{picture}(30,40)
            \qbezier(0,5)(40,20)(0,35)
            \qbezier(30,35)(8,25)(30,35)
            \qbezier(11,26)(1,20)(11,14)
            \qbezier(30,5)(7,14)(30,5)
            {\color{blue}
            {\put(8,27.5){\circle*{3}}
\put(20,27.5){\circle*{3}}
\put(9,27.5){\line(1,0){10}}
            }}
            {\color{red}{\put(16,12){\circle*{3}}
\put(5,12){\circle*{3}}
\put(5,12){\line(1,0){10}}}}
\put(5,17){\text{$+$}}
\put(7,37){$p$}
\put(7,0){$q$}
        \end{picture}
    \end{minipage} \otimes [xb] \mapsto
    - \left(\qquad  \begin{minipage}{40pt}
        \begin{picture}(30,40)
            \qbezier(0,5)(40,20)(0,35)
            \qbezier(30,35)(8,25)(30,35)
            \qbezier(11,26)(1,20)(11,14)
            \qbezier(30,5)(10,15)(30,5)
  {\color{red}{\put(14,33){\circle*{3}}
\put(14,24){\circle*{3}}
\put(14,24){\line(0,1){8}}}}
            {\color{blue}{\put(10,17){\circle*{3}}
\put(10,8){\circle*{3}}
\put(10,9){\line(0,1){8}}}}
\put(-23,19){$p:q$}
\put(20,19){$q:p$}
        \end{picture}
    \end{minipage}\quad \otimes [xa] \quad + \quad 
    \begin{minipage}{40pt}
        \begin{picture}(30,40)
            \qbezier(0,5)(40,20)(0,35)
            \qbezier(30,35)(8,25)(30,35)
            \qbezier(11,26)(1,20)(11,14)
            \qbezier(30,5)(7,14)(30,5)
            {\color{blue}
            {\put(8,27.5){\circle*{3}}
\put(20,27.5){\circle*{3}}
\put(9,27.5){\line(1,0){10}}
            }}
            {\color{red}{\put(16,12){\circle*{3}}
\put(5,12){\circle*{3}}
\put(5,12){\line(1,0){10}}}}
\put(5,17){\text{$-$}}
\put(-25,42){$(p:q):(q:p)$}
\put(-25,-5){$(q:p):(p:q)$}
        \end{picture}
    \end{minipage} \otimes [xb]~\right), \]\[\text{otherwise} \mapsto 0.  \]

\vspace{3mm}
The homotopy $h$ connecting ${\rm in} \circ \rho$ to the identity $:$ $\mathcal{C}\left(~
    \begin{minipage}{30pt}
        \begin{picture}(30,30)
            \qbezier(0,0)(40,15)(0,30)
            \qbezier(30,30)(8,20)(30,30)
            \qbezier(12,20)(2,15)(12,10)
            \qbezier(30,0)(10,10)(30,0)
        \end{picture}
    \end{minipage}
~\right) \to \mathcal{C}\left(~
    \begin{minipage}{30pt}
        \begin{picture}(30,30)
            \qbezier(0,0)(40,15)(0,30)
            \qbezier(30,30)(8,20)(30,30)
            \qbezier(12,20)(2,15)(12,10)
            \qbezier(30,0)(10,10)(30,0)
        \end{picture}
    \end{minipage}
~\right)$ such that $d \circ h$ $+$ $h \circ d$ $=$ $id - in \circ \rho$, is defined by the formulas: 
\begin{align}
&\begin{minipage}{30pt}
        \begin{picture}(30,30)
            \qbezier(0,0)(40,15)(0,30)
            \qbezier(30,30)(8,20)(30,30)
            \qbezier(11,21)(1,15)(11,9)
            \qbezier(30,0)(7,9)(30,0)
  {\color{red}{\put(14,28){\circle*{3}}
\put(14,19){\circle*{3}}
\put(14,19){\line(0,1){8}}}}
            {\color{red}{\put(16,7){\circle*{3}}
\put(5,7){\circle*{3}}
\put(5,7){\line(1,0){10}}}}
\put(8,34){$p$}
\put(8,-3){$q$}
        \end{picture}
    \end{minipage} \otimes [xab] \mapsto - \begin{minipage}{40pt}
        \begin{picture}(30,40)
            \qbezier(0,5)(40,20)(0,35)
            \qbezier(30,35)(8,25)(30,35)
            \qbezier(11,26)(1,20)(11,14)
            \qbezier(30,5)(7,14)(30,5)
            {\color{blue}
            {\put(8,27.5){\circle*{3}}
\put(20,27.5){\circle*{3}}
\put(9,27.5){\line(1,0){10}}
            }}
            {\color{red}{\put(16,12){\circle*{3}}
\put(5,12){\circle*{3}}
\put(5,12){\line(1,0){10}}}}
\put(5,17){\text{$-$}}
\put(7,37){$p$}
\put(7,0){$q$}
        \end{picture}
    \end{minipage} \otimes [xb], \quad  \begin{minipage}{40pt}
        \begin{picture}(30,40)
            \qbezier(0,5)(40,20)(0,35)
            \qbezier(30,35)(8,25)(30,35)
            \qbezier(11,26)(1,20)(11,14)
            \qbezier(30,5)(7,14)(30,5)
            {\color{blue}
            {\put(8,27.5){\circle*{3}}
\put(20,27.5){\circle*{3}}
\put(9,27.5){\line(1,0){10}}
            }}
            {\color{red}{\put(16,12){\circle*{3}}
\put(5,12){\circle*{3}}
\put(5,12){\line(1,0){10}}}}
\put(5,17){\text{$+$}}
\put(7,37){$p$}
\put(7,0){$q$}  
        \end{picture}
    \end{minipage}\otimes [xb] \mapsto \begin{minipage}{30pt}
        \begin{picture}(30,30)
            \qbezier(0,0)(40,15)(0,30)
            \qbezier(30,30)(8,20)(30,30)
            \qbezier(11,21)(1,15)(11,9)
            \qbezier(30,0)(10,10)(30,0)
            {\color{blue}
            {\put(8,22.5){\circle*{3}}
\put(20,22.5){\circle*{3}}
\put(9,22.5){\line(1,0){10}}
            }}
            {\color{blue}{\put(10,12){\circle*{3}}
\put(10,3){\circle*{3}}
\put(10,4){\line(0,1){8}}}}
\put(8,30){$p$}
\put(8,-7){$q$}
        \end{picture}
    \end{minipage} \otimes [x], \\ \nonumber
     &\text{otherwise} \mapsto 0.  
\end{align}

\subsection{Third Reidemeister move}
Let $a$, $b$, $c$ be crossings, $x$ sequence of crossings with negative markers  and $p$, $q$, $r$ be signs.  
For a crossing with no markers or no signs in the following formulas, any markers or signs may be selected.  Let $\tilde{r}$ be $r$ unless the upper left arc is connected to one of the other arcs in the picture.  Let $\tilde{q}$ be $q$ unless the lower right arc is connected to one of the other arcs in the picture.  We use the symbols $p:q$, $q:p$, $p:r$ and $r:p$ whose definitions are given by \cite[Figure 3 Subsection 2.2]{jacobsson}.

\begin{align}
& \mathcal{C}\left(~
    \begin{minipage}{60pt}
        \begin{picture}(50,60)
            \qbezier(0,0)(5,5)(12,12)
            \qbezier(19,19)(25,25)(31,31)
            \qbezier(38,38)(45,45)(50,50)
            \qbezier(0,50)(17.25,44.25)(34.5,34.5)
            \qbezier(34.5,34.5)(44,30)(44,0)
            \qbezier(17,48)(21,53)(24,56)
            \qbezier(12,41)(-1,25)(33.5,0)
            \put(4,40){$c$}
            \put(40,30){$b$}
            \put(3,12){$a$}
        \end{picture}
    \end{minipage}
~\right) = \\ \nonumber
&
\mathcal{C'}\left(~\begin{minipage}{60pt}
        \begin{picture}(50,60)
            \qbezier(0,0)(5,5)(12,12)
            \qbezier(19,19)(25,25)(31,31)
            \qbezier(38,38)(45,45)(50,50)
            \qbezier(0,50)(17.25,44.25)(34.5,34.5)
            \qbezier(34.5,34.5)(44,30)(44,0)
            \qbezier(17,48)(21,53)(24,56)
            \qbezier(12,41)(-1,25)(33.5,0)
            {\color{red}{\put(16,21){\circle*{3}}
\put(16,10){\circle*{3}}
\put(16,11){\line(0,1){10}}}}
{\color{blue}
            {\put(4.5,44.5){\circle*{3}}
\put(16.5,44.5){\circle*{3}}
\put(5.5,44.5){\line(1,0){10}}
            }}
            {\color{blue}
            {\put(32.5,34.5){\circle*{3}}
\put(20.5,34.5){\circle*{3}}
\put(21.5,34.5){\line(1,0){10}}
            }}
            \put(13,14){$q$}
            \put(45,45){$p$}
            \put(4,48){$r$}
        \end{picture}
    \end{minipage} \otimes [xa] + \begin{minipage}{60pt}
        \begin{picture}(50,60)
            \qbezier(0,0)(5,5)(12,12)
            \qbezier(19,19)(25,25)(31,31)
            \qbezier(38,38)(45,45)(50,50)
            \qbezier(0,50)(17.25,44.25)(34.5,34.5)
            \qbezier(34.5,34.5)(44,30)(44,0)
            \qbezier(17,48)(21,53)(24,56)
            \qbezier(12,41)(-1,25)(33.5,0)
            \put(15,28){\text{$-$}}
            {\color{blue}
            {\put(8,44.5){\circle*{3}}
\put(20,44.5){\circle*{3}}
\put(9,44.5){\line(1,0){10}}
            }}
            {\color{blue}{\put(8,15){\circle*{3}}
\put(18,15){\circle*{3}}
\put(9,15){\line(1,0){10}}}}
{\color{red}{\put(27,39){\circle*{3}}
\put(27,29){\circle*{3}}
\put(27,29){\line(0,1){10}}}}
\put(-2,-2){$p:q$}
\put(40,0){$q:p$}
\put(5,50){$\tilde{r}$}
        \end{picture}
    \end{minipage} \otimes [xb], \begin{minipage}{60pt}
        \begin{picture}(50,60)
            \qbezier(0,0)(5,5)(12,12)
            \qbezier(19,19)(25,25)(31,31)
            \qbezier(38,38)(45,45)(50,50)
            \qbezier(0,50)(17.25,44.25)(34.5,34.5)
            \qbezier(34.5,34.5)(44,30)(44,0)
            \qbezier(17,48)(21,53)(24,56)
            \qbezier(12,41)(-1,25)(33.5,0)
            {\color{red}{\put(14,49){\circle*{3}}
\put(14,39){\circle*{3}}
\put(14,39){\line(0,1){10}}}}
        \end{picture}
    \end{minipage} \otimes [x]~\right)\\ \nonumber
    & \qquad \oplus
    \mathcal{C'}_{contr}\left(~ 
    \begin{minipage}{60pt}
        \begin{picture}(50,60)
            \qbezier(0,0)(5,5)(12,12)
            \qbezier(19,19)(25,25)(31,31)
            \qbezier(38,38)(45,45)(50,50)
            \qbezier(0,50)(17.25,44.25)(34.5,34.5)
            \qbezier(34.5,34.5)(44,30)(44,0)
            \qbezier(17,48)(21,53)(24,56)
            \qbezier(12,41)(-1,25)(33.5,0)
            {\color{blue}
            {\put(8,44.5){\circle*{3}}
\put(20,44.5){\circle*{3}}
\put(9,44.5){\line(1,0){10}}
            }}
            {\color{blue}
            {\put(37.5,34.5){\circle*{3}}
\put(25.5,34.5){\circle*{3}}
\put(26.5,34.5){\line(1,0){10}}
            }}
            {\color{blue}{\put(4,15){\circle*{3}}
\put(14,15){\circle*{3}}
\put(5,15){\line(1,0){10}}}}
        \end{picture}
    \end{minipage} \otimes [x] , 
    \begin{minipage}{60pt}
        \begin{picture}(50,60)
        \qbezier(0,0)(5,5)(12,12)
            \qbezier(19,19)(25,25)(31,31)
            \qbezier(38,38)(45,45)(50,50)
            \qbezier(0,50)(17.25,44.25)(34.5,34.5)
            \qbezier(34.5,34.5)(44,30)(44,0)
            \qbezier(17,48)(21,53)(24,56)
            \qbezier(12,41)(-1,25)(33.5,0)
            \put(15,28){\text{$-$}}
            {\color{blue}
            {\put(8,44.5){\circle*{3}}
\put(20,44.5){\circle*{3}}
\put(9,44.5){\line(1,0){10}}
            }}
            {\color{blue}{\put(8,15){\circle*{3}}
\put(18,15){\circle*{3}}
\put(9,15){\line(1,0){10}}}}
{\color{red}{\put(27,39){\circle*{3}}
\put(27,29){\circle*{3}}
\put(27,29){\line(0,1){10}}}}
                    \end{picture}
    \end{minipage} \otimes [xb] ~\right).  \nonumber
\end{align}

\begin{align}
& \mathcal{C}\left(~
    \begin{minipage}{60pt}
        \begin{picture}(50,60)
            \qbezier(50,60)(45,55)(38,48)
            \qbezier(31,41)(25,35)(19,29)
            \qbezier(12,22)(5,15)(0,10)
            \qbezier(50,10)(32.75,15.75)(15.5,25.5)
            \qbezier(15.5,25.5)(6,30)(6,60)
            \qbezier(33,12)(29,7)(26,4)
            \qbezier(38,19)(51,35)(16.5,60)
            \put(26,11){$c$}
            \put(5,23){$a$}
            \put(24,40){$b$}
        \end{picture}
    \end{minipage}
~\right) = \\ \nonumber
&
\mathcal{C}\left(~\begin{minipage}{60pt}
        \begin{picture}(50,60)
            \qbezier(50,60)(45,55)(38,48)
            \qbezier(31,41)(25,35)(19,29)
            \qbezier(12,22)(5,15)(0,10)
            \qbezier(50,10)(32.75,15.75)(15.5,25.5)
            \qbezier(15.5,25.5)(6,30)(6,60)
            \qbezier(33,12)(29,7)(26,4)
            \qbezier(38,19)(51,35)(16.5,60)
            {\color{red}{\put(34,39){\circle*{3}}
\put(34,50){\circle*{3}}
\put(34,49){\line(0,-1){10}}}}
{\color{blue}
            {\put(38.5,15.5){\circle*{3}}
\put(27.5,15.5){\circle*{3}}
\put(37.5,15.5){\line(-1,0){10}}
            }}
            {\color{blue}
            {\put(2.5,25.5){\circle*{3}}
\put(14.5,25.5){\circle*{3}}
\put(13.5,25.5){\line(-1,0){10}}
            }}
            \put(27,5){$q$}
            \put(35,45){$p$}
            \put(0,53){$r$}
        \end{picture}
    \end{minipage} \otimes [xb] + \begin{minipage}{60pt}
        \begin{picture}(60,60)
            \qbezier(50,60)(45,55)(38,48)
            \qbezier(31,41)(25,35)(19,29)
            \qbezier(12,22)(5,15)(0,10)
            \qbezier(50,10)(32.75,15.75)(15.5,25.5)
            \qbezier(15.5,25.5)(6,30)(6,60)
            \qbezier(33,12)(29,7)(26,4)
            \qbezier(38,19)(51,35)(16.5,60)
            \put(27,25){\text{$-$}}
            {\color{blue}
            {\put(42,15.5){\circle*{3}}
\put(30,15.5){\circle*{3}}
\put(41,15.5){\line(-1,0){10}}
            }}
            {\color{blue}{\put(34,45){\circle*{3}}
\put(24,45){\circle*{3}}
\put(33,45){\line(-1,0){10}}}}
{\color{red}{\put(7,21){\circle*{3}}
\put(7,31){\circle*{3}}
\put(7,31){\line(0,-1){10}}}}
\put(17,60){$p:r$}
\put(-25,57){$r:p$}
\put(25,2){$\tilde{q}$}
        \end{picture}
    \end{minipage} \otimes [xa], \begin{minipage}{60pt}
        \begin{picture}(50,60)
            \qbezier(50,60)(45,55)(38,48)
            \qbezier(31,41)(25,35)(19,29)
            \qbezier(12,22)(5,15)(0,10)
            \qbezier(50,10)(32.75,15.75)(15.5,25.5)
            \qbezier(15.5,25.5)(6,30)(6,60)
            \qbezier(33,12)(29,7)(26,4)
            \qbezier(38,19)(51,35)(16.5,60)
            {\color{red}{\put(36,11){\circle*{3}}
\put(36,21){\circle*{3}}
\put(36,21){\line(0,-1){10}}}}
        \end{picture}
    \end{minipage} \otimes [x]~\right)\\ \nonumber
    & \qquad \oplus
    \mathcal{C}_{contr}\left(~ 
    \begin{minipage}{60pt}
        \begin{picture}(50,60)
            \qbezier(50,60)(45,55)(38,48)
            \qbezier(31,41)(25,35)(19,29)
            \qbezier(12,22)(5,15)(0,10)
            \qbezier(50,10)(32.75,15.75)(15.5,25.5)
            \qbezier(15.5,25.5)(6,30)(6,60)
            \qbezier(33,12)(29,7)(26,4)
            \qbezier(38,19)(51,35)(16.5,60)
            {\color{blue}
            {\put(42,15.5){\circle*{3}}
\put(30,15.5){\circle*{3}}
\put(41,15.5){\line(-1,0){10}}
            }}
            {\color{blue}
            {\put(7,25.5){\circle*{3}}
\put(17,25.5){\circle*{3}}
\put(17,25.5){\line(-1,0){10}}
            }}
            {\color{blue}{\put(31,45){\circle*{3}}
\put(21,45){\circle*{3}}
\put(30,45){\line(-1,0){10}}}}
        \end{picture}
    \end{minipage} \otimes [x] , 
    \begin{minipage}{60pt}
        \begin{picture}(50,60)
        \qbezier(50,60)(45,55)(38,48)
            \qbezier(31,41)(25,35)(19,29)
            \qbezier(12,22)(5,15)(0,10)
            \qbezier(50,10)(32.75,15.75)(15.5,25.5)
            \qbezier(15.5,25.5)(6,30)(6,60)
            \qbezier(33,12)(29,7)(26,4)
            \qbezier(38,19)(51,35)(16.5,60)
            \put(27,25){\text{$-$}}
            {\color{blue}
            {\put(42,15.5){\circle*{3}}
\put(30,15.5){\circle*{3}}
\put(41,15.5){\line(-1,0){10}}
            }}
            {\color{blue}{\put(34,45){\circle*{3}}
\put(24,45){\circle*{3}}
\put(33,45){\line(-1,0){10}}}}
{\color{red}{\put(7,21){\circle*{3}}
\put(7,31){\circle*{3}}
\put(7,31){\line(0,-1){10}}}}
                    \end{picture}
    \end{minipage} \otimes [xa] ~\right).  \nonumber
\end{align}

Let a link diagram $D'$ $=$ \begin{minipage}{60pt}
        \begin{picture}(50,60)
            \qbezier(50,60)(45,55)(38,48)
            \qbezier(31,41)(25,35)(19,29)
            \qbezier(12,22)(5,15)(0,10)
            \qbezier(50,10)(32.75,15.75)(15.5,25.5)
            \qbezier(15.5,25.5)(6,30)(6,60)
            \qbezier(33,12)(29,7)(26,4)
            \qbezier(38,19)(51,35)(16.5,60)
            \put(26,11){$c$}
            \put(5,23){$a$}
            \put(24,40){$b$}
        \end{picture}
    \end{minipage} and $D$ $=$
    \begin{minipage}{60pt}
        \begin{picture}(50,60)
            \qbezier(50,60)(45,55)(38,48)
            \qbezier(31,41)(25,35)(19,29)
            \qbezier(12,22)(5,15)(0,10)
            \qbezier(50,10)(32.75,15.75)(15.5,25.5)
            \qbezier(15.5,25.5)(6,30)(6,60)
            \qbezier(33,12)(29,7)(26,4)
            \qbezier(38,19)(51,35)(16.5,60)
            \put(26,11){$c$}
            \put(5,23){$a$}
            \put(24,40){$b$}
        \end{picture}
    \end{minipage}.  

By using above formula, consider the following mapping 
\[C(D') = C' \oplus C'_{contr} \stackrel{\rho}{\to} C' \stackrel{{\rm isom}}{\to} C \stackrel{i}{\to} C \oplus C_{contr} = C(D) \]
as in \cite[Page 1223]{jacobsson}.  

The isomorphism $C' \to C$ is defined formulas

\begin{align}
&\begin{minipage}{60pt}
        \begin{picture}(50,60)
            \qbezier(0,0)(5,5)(12,12)
            \qbezier(19,19)(25,25)(31,31)
            \qbezier(38,38)(45,45)(50,50)
            \qbezier(0,50)(17.25,44.25)(34.5,34.5)
            \qbezier(34.5,34.5)(44,30)(44,0)
            \qbezier(17,48)(21,53)(24,56)
            \qbezier(12,41)(-1,25)(33.5,0)
            {\color{red}{\put(16,21){\circle*{3}}
\put(16,10){\circle*{3}}
\put(16,11){\line(0,1){10}}}}
{\color{blue}
            {\put(4.5,44.5){\circle*{3}}
\put(16.5,44.5){\circle*{3}}
\put(5.5,44.5){\line(1,0){10}}
            }}
            {\color{blue}
            {\put(32.5,34.5){\circle*{3}}
\put(20.5,34.5){\circle*{3}}
\put(21.5,34.5){\line(1,0){10}}
            }}
            \put(13,14){$q$}
            \put(45,45){$p$}
            \put(4,48){$r$}
        \end{picture}
    \end{minipage} \otimes [xa] + \begin{minipage}{60pt}
        \begin{picture}(50,60)
            \qbezier(0,0)(5,5)(12,12)
            \qbezier(19,19)(25,25)(31,31)
            \qbezier(38,38)(45,45)(50,50)
            \qbezier(0,50)(17.25,44.25)(34.5,34.5)
            \qbezier(34.5,34.5)(44,30)(44,0)
            \qbezier(17,48)(21,53)(24,56)
            \qbezier(12,41)(-1,25)(33.5,0)
            \put(15,28){\text{$-$}}
            {\color{blue}
            {\put(8,44.5){\circle*{3}}
\put(20,44.5){\circle*{3}}
\put(9,44.5){\line(1,0){10}}
            }}
            {\color{blue}{\put(8,15){\circle*{3}}
\put(18,15){\circle*{3}}
\put(9,15){\line(1,0){10}}}}
{\color{red}{\put(27,39){\circle*{3}}
\put(27,29){\circle*{3}}
\put(27,29){\line(0,1){10}}}}
\put(-2,-2){$p:q$}
\put(40,0){$q:p$}
\put(5,50){$\tilde{r}$}
        \end{picture}
    \end{minipage} \otimes [xb] \mapsto \begin{minipage}{60pt}
        \begin{picture}(50,60)
            \qbezier(50,60)(45,55)(38,48)
            \qbezier(31,41)(25,35)(19,29)
            \qbezier(12,22)(5,15)(0,10)
            \qbezier(50,10)(32.75,15.75)(15.5,25.5)
            \qbezier(15.5,25.5)(6,30)(6,60)
            \qbezier(33,12)(29,7)(26,4)
            \qbezier(38,19)(51,35)(16.5,60)
            {\color{red}{\put(34,39){\circle*{3}}
\put(34,50){\circle*{3}}
\put(34,49){\line(0,-1){10}}}}
{\color{blue}
            {\put(38.5,15.5){\circle*{3}}
\put(27.5,15.5){\circle*{3}}
\put(37.5,15.5){\line(-1,0){10}}
            }}
            {\color{blue}
            {\put(2.5,25.5){\circle*{3}}
\put(14.5,25.5){\circle*{3}}
\put(13.5,25.5){\line(-1,0){10}}
            }}
            \put(27,5){$q$}
            \put(35,45){$p$}
            \put(0,53){$r$}
        \end{picture}
    \end{minipage} \otimes [xb] + \begin{minipage}{60pt}
        \begin{picture}(60,60)
            \qbezier(50,60)(45,55)(38,48)
            \qbezier(31,41)(25,35)(19,29)
            \qbezier(12,22)(5,15)(0,10)
            \qbezier(50,10)(32.75,15.75)(15.5,25.5)
            \qbezier(15.5,25.5)(6,30)(6,60)
            \qbezier(33,12)(29,7)(26,4)
            \qbezier(38,19)(51,35)(16.5,60)
            \put(27,25){\text{$-$}}
            {\color{blue}
            {\put(42,15.5){\circle*{3}}
\put(30,15.5){\circle*{3}}
\put(41,15.5){\line(-1,0){10}}
            }}
            {\color{blue}{\put(34,45){\circle*{3}}
\put(24,45){\circle*{3}}
\put(33,45){\line(-1,0){10}}}}
{\color{red}{\put(7,21){\circle*{3}}
\put(7,31){\circle*{3}}
\put(7,31){\line(0,-1){10}}}}
\put(17,60){$p:r$}
\put(-25,57){$r:p$}
\put(25,2){$\tilde{q}$}
        \end{picture}
    \end{minipage} \otimes [xa], \\ \nonumber 
    \\ \nonumber
    &\begin{minipage}{60pt}
        \begin{picture}(50,60)
            \qbezier(0,0)(5,5)(12,12)
            \qbezier(19,19)(25,25)(31,31)
            \qbezier(38,38)(45,45)(50,50)
            \qbezier(0,50)(17.25,44.25)(34.5,34.5)
            \qbezier(34.5,34.5)(44,30)(44,0)
            \qbezier(17,48)(21,53)(24,56)
            \qbezier(12,41)(-1,25)(33.5,0)
            {\color{red}{\put(14,49){\circle*{3}}
\put(14,39){\circle*{3}}
\put(14,39){\line(0,1){10}}}}
        \end{picture}
    \end{minipage} \otimes [x] \mapsto \begin{minipage}{60pt}
        \begin{picture}(50,60)
            \qbezier(50,60)(45,55)(38,48)
            \qbezier(31,41)(25,35)(19,29)
            \qbezier(12,22)(5,15)(0,10)
            \qbezier(50,10)(32.75,15.75)(15.5,25.5)
            \qbezier(15.5,25.5)(6,30)(6,60)
            \qbezier(33,12)(29,7)(26,4)
            \qbezier(38,19)(51,35)(16.5,60)
            {\color{red}{\put(36,11){\circle*{3}}
\put(36,21){\circle*{3}}
\put(36,21){\line(0,-1){10}}}}
        \end{picture}
    \end{minipage} \otimes [x].  
\end{align}

The retraction $\rho : \mathcal{C}(D') \to
\mathcal{C'}\left(~\begin{minipage}{60pt}
        \begin{picture}(50,60)
            \qbezier(0,0)(5,5)(12,12)
            \qbezier(19,19)(25,25)(31,31)
            \qbezier(38,38)(45,45)(50,50)
            \qbezier(0,50)(17.25,44.25)(34.5,34.5)
            \qbezier(34.5,34.5)(44,30)(44,0)
            \qbezier(17,48)(21,53)(24,56)
            \qbezier(12,41)(-1,25)(33.5,0)
            {\color{red}{\put(16,21){\circle*{3}}
\put(16,10){\circle*{3}}
\put(16,11){\line(0,1){10}}}}
{\color{blue}
            {\put(4.5,44.5){\circle*{3}}
\put(16.5,44.5){\circle*{3}}
\put(5.5,44.5){\line(1,0){10}}
            }}
            {\color{blue}
            {\put(32.5,34.5){\circle*{3}}
\put(20.5,34.5){\circle*{3}}
\put(21.5,34.5){\line(1,0){10}}
            }}
            \put(13,14){$q$}
            \put(45,45){$p$}
            \put(4,48){$r$}
        \end{picture}
    \end{minipage} \otimes [xa] + \begin{minipage}{60pt}
        \begin{picture}(50,60)
            \qbezier(0,0)(5,5)(12,12)
            \qbezier(19,19)(25,25)(31,31)
            \qbezier(38,38)(45,45)(50,50)
            \qbezier(0,50)(17.25,44.25)(34.5,34.5)
            \qbezier(34.5,34.5)(44,30)(44,0)
            \qbezier(17,48)(21,53)(24,56)
            \qbezier(12,41)(-1,25)(33.5,0)
            \put(15,28){\text{$-$}}
            {\color{blue}
            {\put(8,44.5){\circle*{3}}
\put(20,44.5){\circle*{3}}
\put(9,44.5){\line(1,0){10}}
            }}
            {\color{blue}{\put(8,15){\circle*{3}}
\put(18,15){\circle*{3}}
\put(9,15){\line(1,0){10}}}}
{\color{red}{\put(27,39){\circle*{3}}
\put(27,29){\circle*{3}}
\put(27,29){\line(0,1){10}}}}
\put(-2,-2){$p:q$}
\put(40,0){$q:p$}
\put(5,50){$\tilde{r}$}
        \end{picture}
    \end{minipage} \otimes [xb], \begin{minipage}{60pt}
        \begin{picture}(50,60)
            \qbezier(0,0)(5,5)(12,12)
            \qbezier(19,19)(25,25)(31,31)
            \qbezier(38,38)(45,45)(50,50)
            \qbezier(0,50)(17.25,44.25)(34.5,34.5)
            \qbezier(34.5,34.5)(44,30)(44,0)
            \qbezier(17,48)(21,53)(24,56)
            \qbezier(12,41)(-1,25)(33.5,0)
            {\color{red}{\put(14,49){\circle*{3}}
\put(14,39){\circle*{3}}
\put(14,39){\line(0,1){10}}}}
        \end{picture}
    \end{minipage} \otimes [x]~\right)$ is defined by the formulas 
\begin{align}
&\begin{minipage}{60pt}
        \begin{picture}(50,60)
            \qbezier(0,0)(5,5)(12,12)
            \qbezier(19,19)(25,25)(31,31)
            \qbezier(38,38)(45,45)(50,50)
            \qbezier(0,50)(17.25,44.25)(34.5,34.5)
            \qbezier(34.5,34.5)(44,30)(44,0)
            \qbezier(17,48)(21,53)(24,56)
            \qbezier(12,41)(-1,25)(33.5,0)
            {\color{red}{\put(16,21){\circle*{3}}
\put(16,10){\circle*{3}}
\put(16,11){\line(0,1){10}}}}
{\color{blue}
            {\put(4.5,44.5){\circle*{3}}
\put(16.5,44.5){\circle*{3}}
\put(5.5,44.5){\line(1,0){10}}
            }}
            {\color{blue}
            {\put(32.5,34.5){\circle*{3}}
\put(20.5,34.5){\circle*{3}}
\put(21.5,34.5){\line(1,0){10}}
            }}
            \put(13,14){$q$}
            \put(45,45){$p$}
            \put(4,48){$r$}
        \end{picture}
    \end{minipage} \otimes [xa] \mapsto \begin{minipage}{60pt}
        \begin{picture}(50,60)
            \qbezier(0,0)(5,5)(12,12)
            \qbezier(19,19)(25,25)(31,31)
            \qbezier(38,38)(45,45)(50,50)
            \qbezier(0,50)(17.25,44.25)(34.5,34.5)
            \qbezier(34.5,34.5)(44,30)(44,0)
            \qbezier(17,48)(21,53)(24,56)
            \qbezier(12,41)(-1,25)(33.5,0)
            {\color{red}{\put(16,21){\circle*{3}}
\put(16,10){\circle*{3}}
\put(16,11){\line(0,1){10}}}}
{\color{blue}
            {\put(4.5,44.5){\circle*{3}}
\put(16.5,44.5){\circle*{3}}
\put(5.5,44.5){\line(1,0){10}}
            }}
            {\color{blue}
            {\put(32.5,34.5){\circle*{3}}
\put(20.5,34.5){\circle*{3}}
\put(21.5,34.5){\line(1,0){10}}
            }}
            \put(13,14){$q$}
            \put(45,45){$p$}
            \put(4,48){$r$}
        \end{picture}
    \end{minipage} \otimes [xa] + \begin{minipage}{60pt}
        \begin{picture}(50,60)
            \qbezier(0,0)(5,5)(12,12)
            \qbezier(19,19)(25,25)(31,31)
            \qbezier(38,38)(45,45)(50,50)
            \qbezier(0,50)(17.25,44.25)(34.5,34.5)
            \qbezier(34.5,34.5)(44,30)(44,0)
            \qbezier(17,48)(21,53)(24,56)
            \qbezier(12,41)(-1,25)(33.5,0)
            \put(15,28){\text{$-$}}
            {\color{blue}
            {\put(8,44.5){\circle*{3}}
\put(20,44.5){\circle*{3}}
\put(9,44.5){\line(1,0){10}}
            }}
            {\color{blue}{\put(8,15){\circle*{3}}
\put(18,15){\circle*{3}}
\put(9,15){\line(1,0){10}}}}
{\color{red}{\put(27,39){\circle*{3}}
\put(27,29){\circle*{3}}
\put(27,29){\line(0,1){10}}}}
\put(-2,-2){$p:q$}
\put(40,0){$q:p$}
\put(5,50){$\tilde{r}$}
        \end{picture}
    \end{minipage} \otimes [xb], \\ \nonumber
&\begin{minipage}{60pt}
        \begin{picture}(50,60)
            \qbezier(0,0)(5,5)(12,12)
            \qbezier(19,19)(25,25)(31,31)
            \qbezier(38,38)(45,45)(50,50)
            \qbezier(0,50)(17.25,44.25)(34.5,34.5)
            \qbezier(34.5,34.5)(44,30)(44,0)
            \qbezier(17,48)(21,53)(24,56)
            \qbezier(12,41)(-1,25)(33.5,0)
            {\color{red}{\put(14,49){\circle*{3}}
\put(14,39){\circle*{3}}
\put(14,39){\line(0,1){10}}}}
        \end{picture}
    \end{minipage} \otimes [x] \mapsto \begin{minipage}{60pt}
        \begin{picture}(50,60)
            \qbezier(0,0)(5,5)(12,12)
            \qbezier(19,19)(25,25)(31,31)
            \qbezier(38,38)(45,45)(50,50)
            \qbezier(0,50)(17.25,44.25)(34.5,34.5)
            \qbezier(34.5,34.5)(44,30)(44,0)
            \qbezier(17,48)(21,53)(24,56)
            \qbezier(12,41)(-1,25)(33.5,0)
            {\color{red}{\put(14,49){\circle*{3}}
\put(14,39){\circle*{3}}
\put(14,39){\line(0,1){10}}}}
        \end{picture}
    \end{minipage} \otimes [x], \\ \nonumber
\\ \nonumber   
&\begin{minipage}{60pt}
        \begin{picture}(50,60)
            \qbezier(0,0)(5,5)(12,12)
            \qbezier(19,19)(25,25)(31,31)
            \qbezier(38,38)(45,45)(50,50)
            \qbezier(0,50)(17.25,44.25)(34.5,34.5)
            \qbezier(34.5,34.5)(44,30)(44,0)
            \qbezier(17,48)(21,53)(24,56)
            \qbezier(12,41)(-1,25)(33.5,0)
            \put(15,28){\text{$+$}}
            {\color{blue}
            {\put(8,44.5){\circle*{3}}
\put(20,44.5){\circle*{3}}
\put(9,44.5){\line(1,0){10}}
            }}
            {\color{blue}{\put(8,15){\circle*{3}}
\put(18,15){\circle*{3}}
\put(9,15){\line(1,0){10}}}}
{\color{red}{\put(27,39){\circle*{3}}
\put(27,29){\circle*{3}}
\put(27,29){\line(0,1){10}}}}
\put(5,-2){$q$}
\put(40,0){$p$}
\put(5,50){$r$}
        \end{picture}
    \end{minipage}\!\!\!\!\! \otimes [xb] \mapsto - \begin{minipage}{60pt}
        \begin{picture}(50,60)
            \qbezier(0,0)(5,5)(12,12)
            \qbezier(19,19)(25,25)(31,31)
            \qbezier(38,38)(45,45)(50,50)
            \qbezier(0,50)(17.25,44.25)(34.5,34.5)
            \qbezier(34.5,34.5)(44,30)(44,0)
            \qbezier(17,48)(21,53)(24,56)
            \qbezier(12,41)(-1,25)(33.5,0)
            {\color{red}{\put(16,21){\circle*{3}}
\put(16,10){\circle*{3}}
\put(16,11){\line(0,1){10}}}}
{\color{blue}
            {\put(4.5,44.5){\circle*{3}}
\put(16.5,44.5){\circle*{3}}
\put(5.5,44.5){\line(1,0){10}}
            }}
            {\color{blue}
            {\put(32.5,34.5){\circle*{3}}
\put(20.5,34.5){\circle*{3}}
\put(21.5,34.5){\line(1,0){10}}
            }}
            \put(13,14){$p:q$}
            \put(45,45){$q:p$}
            \put(4,48){$r$}
        \end{picture}
    \end{minipage} \otimes [xa] - \begin{minipage}{60pt}
        \begin{picture}(50,60)
            \qbezier(0,0)(5,5)(12,12)
            \qbezier(19,19)(25,25)(31,31)
            \qbezier(38,38)(45,45)(50,50)
            \qbezier(0,50)(17.25,44.25)(34.5,34.5)
            \qbezier(34.5,34.5)(44,30)(44,0)
            \qbezier(17,48)(21,53)(24,56)
            \qbezier(12,41)(-1,25)(33.5,0)
            \put(15,28){\text{$-$}}
            {\color{blue}
            {\put(8,44.5){\circle*{3}}
\put(20,44.5){\circle*{3}}
\put(9,44.5){\line(1,0){10}}
            }}
            {\color{blue}{\put(8,15){\circle*{3}}
\put(18,15){\circle*{3}}
\put(9,15){\line(1,0){10}}}}
{\color{red}{\put(27,39){\circle*{3}}
\put(27,29){\circle*{3}}
\put(27,29){\line(0,1){10}}}}
\put(-32,-13){$(q:p):(p:q)$}
\put(30,55){$(p:q):(q:p)$}
\put(5,50){$\tilde{r}$}
        \end{picture}
    \end{minipage} \otimes [xb] \quad - \begin{minipage}{60pt}
        \begin{picture}(50,60)
            \qbezier(0,0)(5,5)(12,12)
            \qbezier(19,19)(25,25)(31,31)
            \qbezier(38,38)(45,45)(50,50)
            \qbezier(0,50)(17.25,44.25)(34.5,34.5)
            \qbezier(34.5,34.5)(44,30)(44,0)
            \qbezier(17,48)(21,53)(24,56)
            \qbezier(12,41)(-1,25)(33.5,0)
            {\color{red}
            {\put(15,50){\circle*{3}}
\put(15,40){\circle*{3}}
\put(15,40){\line(0,1){10}}
            }}
            {\color{blue}{\put(8,15){\circle*{3}}
\put(18,15){\circle*{3}}
\put(9,15){\line(1,0){10}}}}
{\color{blue}{\put(32,34){\circle*{3}}
\put(22,34){\circle*{3}}
\put(22,34){\line(1,0){10}}}}
\put(5,-2){$q$}
\put(40,55){$p:r$}
\put(-20,40){$r:p$}
        \end{picture}
    \end{minipage} \otimes [xc], \\ \nonumber
\\ \nonumber    
&\begin{minipage}{60pt}
        \begin{picture}(50,60)
            \qbezier(0,0)(5,5)(12,12)
            \qbezier(19,19)(25,25)(31,31)
            \qbezier(38,38)(45,45)(50,50)
            \qbezier(0,50)(17.25,44.25)(34.5,34.5)
            \qbezier(34.5,34.5)(44,30)(44,0)
            \qbezier(17,48)(21,53)(24,56)
            \qbezier(12,41)(-1,25)(33.5,0)
            {\color{red}{\put(16,21){\circle*{3}}
\put(16,10){\circle*{3}}
\put(16,11){\line(0,1){10}}}}
{\color{blue}
            {\put(4.5,44.5){\circle*{3}}
\put(16.5,44.5){\circle*{3}}
\put(5.5,44.5){\line(1,0){10}}
            }}
            {\color{red}{\put(27,39){\circle*{3}}
\put(27,29){\circle*{3}}
\put(27,29){\line(0,1){10}}}}
            \put(-3,14){$q$}
            \put(45,45){$p$}
            \put(4,48){$r$}
        \end{picture}
    \end{minipage} \otimes [xab] \mapsto \begin{minipage}{60pt}
        \begin{picture}(50,60)
            \qbezier(0,0)(5,5)(12,12)
            \qbezier(19,19)(25,25)(31,31)
            \qbezier(38,38)(45,45)(50,50)
            \qbezier(0,50)(17.25,44.25)(34.5,34.5)
            \qbezier(34.5,34.5)(44,30)(44,0)
            \qbezier(17,48)(21,53)(24,56)
            \qbezier(12,41)(-1,25)(33.5,0)
            {\color{red}
            {\put(15,50){\circle*{3}}
\put(15,40){\circle*{3}}
\put(15,40){\line(0,1){10}}
            }}
            {\color{blue}{\put(8,15){\circle*{3}}
\put(18,15){\circle*{3}}
\put(9,15){\line(1,0){10}}}}
{\color{red}{\put(27,39){\circle*{3}}
\put(27,29){\circle*{3}}
\put(27,29){\line(0,1){10}}}}
\put(5,-2){$q$}
\put(45,45){$p$}
            \put(-3,40){$r$}
        \end{picture}
    \end{minipage} \otimes [xbc], \\ \nonumber
\\ \nonumber     
&\text{otherwise} \mapsto 0.  
\end{align}
\vspace{3mm}

The homotopy connecting ${\rm in} \circ {\rm isom}^{-1} \circ i^{-1} \circ i \circ {\rm isom} \circ \rho$ ($=$ ${\rm in} \circ \rho$) to identity, that is, a map $h : $ $C(D') \to C(D')$ such that $d \circ h$ $+$ $h \circ d$ $=$ ${\rm id} - {\rm in} \circ \rho$, is defined by the formulas:  
\begin{align}
&\begin{minipage}{60pt}
        \begin{picture}(50,60)
            \qbezier(0,0)(5,5)(12,12)
            \qbezier(19,19)(25,25)(31,31)
            \qbezier(38,38)(45,45)(50,50)
            \qbezier(0,50)(17.25,44.25)(34.5,34.5)
            \qbezier(34.5,34.5)(44,30)(44,0)
            \qbezier(17,48)(21,53)(24,56)
            \qbezier(12,41)(-1,25)(33.5,0)
            {\color{red}{\put(16,21){\circle*{3}}
\put(16,10){\circle*{3}}
\put(16,11){\line(0,1){10}}}}
{\color{blue}
            {\put(4.5,44.5){\circle*{3}}
\put(16.5,44.5){\circle*{3}}
\put(5.5,44.5){\line(1,0){10}}
            }}
            {\color{red}{\put(27,39){\circle*{3}}
\put(27,29){\circle*{3}}
\put(27,29){\line(0,1){10}}}}
            \put(-3,14){$q$}
            \put(45,45){$p$}
            \put(4,48){$r$}
        \end{picture}
    \end{minipage} \otimes [xab] \mapsto -  \begin{minipage}{60pt}
        \begin{picture}(50,60)
            \qbezier(0,0)(5,5)(12,12)
            \qbezier(19,19)(25,25)(31,31)
            \qbezier(38,38)(45,45)(50,50)
            \qbezier(0,50)(17.25,44.25)(34.5,34.5)
            \qbezier(34.5,34.5)(44,30)(44,0)
            \qbezier(17,48)(21,53)(24,56)
            \qbezier(12,41)(-1,25)(33.5,0)
            \put(15,28){\text{$-$}}
            {\color{blue}
            {\put(8,44.5){\circle*{3}}
\put(20,44.5){\circle*{3}}
\put(9,44.5){\line(1,0){10}}
            }}
            {\color{blue}{\put(8,15){\circle*{3}}
\put(18,15){\circle*{3}}
\put(9,15){\line(1,0){10}}}}
{\color{red}{\put(27,39){\circle*{3}}
\put(27,29){\circle*{3}}
\put(27,29){\line(0,1){10}}}}
\put(5,-2){$q$}
            \put(45,45){$p$}
\put(5,50){$r$}
        \end{picture}
    \end{minipage} \otimes [xb], 
    \begin{minipage}{60pt}
        \begin{picture}(50,60)
            \qbezier(0,0)(5,5)(12,12)
            \qbezier(19,19)(25,25)(31,31)
            \qbezier(38,38)(45,45)(50,50)
            \qbezier(0,50)(17.25,44.25)(34.5,34.5)
            \qbezier(34.5,34.5)(44,30)(44,0)
            \qbezier(17,48)(21,53)(24,56)
            \qbezier(12,41)(-1,25)(33.5,0)
            \put(15,28){\text{$+$}}
            {\color{blue}
            {\put(8,44.5){\circle*{3}}
\put(20,44.5){\circle*{3}}
\put(9,44.5){\line(1,0){10}}
            }}
            {\color{blue}{\put(8,15){\circle*{3}}
\put(18,15){\circle*{3}}
\put(9,15){\line(1,0){10}}}}
{\color{red}{\put(27,39){\circle*{3}}
\put(27,29){\circle*{3}}
\put(27,29){\line(0,1){10}}}}
\put(5,-2){$q$}
\put(40,0){$p$}
\put(5,50){$r$}
        \end{picture}
    \end{minipage}\!\!\!\!\!\! \otimes [xb] \mapsto \begin{minipage}{60pt}
        \begin{picture}(50,60)
            \qbezier(0,0)(5,5)(12,12)
            \qbezier(19,19)(25,25)(31,31)
            \qbezier(38,38)(45,45)(50,50)
            \qbezier(0,50)(17.25,44.25)(34.5,34.5)
            \qbezier(34.5,34.5)(44,30)(44,0)
            \qbezier(17,48)(21,53)(24,56)
            \qbezier(12,41)(-1,25)(33.5,0)
            {\color{blue}
            {\put(8,44.5){\circle*{3}}
\put(20,44.5){\circle*{3}}
\put(9,44.5){\line(1,0){10}}
            }}
            {\color{blue}
            {\put(37.5,34.5){\circle*{3}}
\put(25.5,34.5){\circle*{3}}
\put(26.5,34.5){\line(1,0){10}}
            }}
            {\color{blue}{\put(4,15){\circle*{3}}
\put(14,15){\circle*{3}}
\put(5,15){\line(1,0){10}}}}
\put(5,-2){$q$}
\put(40,30){$p$}
\put(5,50){$r$}
        \end{picture}
    \end{minipage} \otimes [x], \\ \nonumber
    \\ \nonumber
&\text{otherwise} \mapsto 0 
\end{align}

\vspace{0.3cm}

Department of Pure and Applied Mathematics, Waseda University

Tokyo 169-8555, JAPAN

email: noboru@moegi.waseda.jp

\end{document}